%
%
%
\input amstex
\documentstyle{amsppt}
\def\twoone{${\hbox{\uppercase\expandafter
{\romannumeral2}}}_1$}
\def\bt{\mathop\boxtimes\limits}
\def\eps{\epsilon}

\def\Aut{{\hbox{Aut}}}
\def\Out{{\hbox{\text{\rm Out}}}}

\def\Int{{\hbox{\rm Int}}}

\def\tr{\hbox{\rm tr}}

\def\op{{\text{\rm op}}}

\TagsOnRight

\baselineskip=18 pt plus 2pt
\lineskip=9 pt minus 1 pt
\lineskiplimit=9 pt

\magnification = \magstep 1


\topmatter
\title
EXAMPLES OF SUBFACTORS WITH PROPERTY T STANDARD INVARIANT
\endtitle
\rightheadtext{Subfactors with property T}
\author DIETMAR BISCH$^1$ AND SORIN POPA$^2$
\endauthor
\leftheadtext{Dietmar Bisch, Sorin Popa}
\thanks $^1$ supported by NSF grant DMS-9531566,
$^2$ supported by NSF grant DMS-9500882 
\endthanks
\affil UC Santa Barbara \\
Department of Mathematics \\
Santa Barbara, CA 93106 \\
  \\
Universit{\'e} de Gen{\`e}ve $\quad \qquad$ \hskip 1in  \ \ UCLA 
$\qquad \qquad$\\
Section de math\'ematiques \hskip 0.3in and \hskip 0.3in 
Department of Mathematics \\
CH-1211 Gen\`eve 24 \ \ \  \hskip 1in \ \ \  Los Angeles, CA 90024
\endaffil
\address \email bisch\@math.ucsb.edu, Sorin.Popa\@math.unige.ch ,
popa\@math.ucla.edu
\endemail \endaddress

\address UC Santa Barbara, Department of Mathematics, Santa Barbara, 
CA 93106, USA and
Universit\'e de Gen\`eve, Section de math\'ematiques, 2-4, rue du Lievre,
Case postale 124, CH-1211 Gen\`eve 24, Switzerland
\endaddress

\subjclass 46L10, 46L37 \endsubjclass

\abstract
Let $H$ and $K$ be two finite groups with a properly outer action
on the II$_1$ factor $M$. We prove that the group type inclusions 
$M^H \subset M \rtimes K$, studied in detail in [BH],
have property $T$ in the sense of [Po6] if and only if
the group generated by $H$ and $K$ in the outer automorphism
group of $M$ has Kazhdan's property T [K]. This construction
yields irreducible, infinite depth subfactors with small Jones indices  
and property T standard invariant.
\endabstract

\endtopmatter

\document

\heading
0. Introduction
\endheading

If $H$ and $K$ are two finite groups with a properly outer action
on the II$_1$ factor $M$, we can compose the two subfactors
$M^H \subset M$ and $M \subset M \rtimes K$ to obtain a new
inclusion $M^H \subset M \rtimes K$. 
While the Jones index of $M^H \subset M \rtimes K$ is finite,
being equal to $|H| \cdot |K|$, in general this inclusion can no longer
be obtained as a fixed point algebra (or a crossed product) by
a group action, or more generally a coaction of a Kac algebra.
It was shown in [BH] that several interesting properties of these group 
type inclusions are determined by the group $G$ generated by
$H$ and $K$ in the outer automorphism group of the factor $M$.
For instance, it was shown there that $M^H \subset M \rtimes K$
is amenable (resp. has finite depth) if and only if $G$ is an amenable
(resp. finite) group. Since any second countable discrete group has
a properly outer action on the hyperfinite II$_1$ factor $R$, numerous
examples of subfactors, whose standard
invariants are essentially as badly or as well behaved as the group $G$
itself, can be constructed in this way.

The second author introduced recently a concept of property T
for the standard invariant of a subfactor and used certain locally
trivial, diagonal subfactors associated to a $G$-kernel to construct
reducible, infinite depth subfactors with finite Jones index whose
standard invariant has property T ([Po5], [Po6]). We provide
in this note new examples of property T subfactors by showing
that the group type inclusions $M^H \subset M \rtimes K$ have
property T in the sense of [Po6] if and only if the
group $G$ has Kazhdan's property T. Once this result is established,
every infinite property $T$ group, which is a quotient of a free
product of two finite groups, provides an example of an {\it irreducible,
infinite depth} subfactor of the hyperfinite II$_1$ factor with
finite Jones index whose standard invariant has property T. The lowest
possible index of an irreducible, infinite depth subfactor with 
property T standard invariant that this
construction can give is index $6$ and we give an example of such
a subfactor in section 2, based on work by Conder.

Here is a more detailed description of the two sections below. 
In section 1 we discuss briefly the different notions of property T
used in this paper and recall a few facts that are needed in the
subsequent section. In section 2 we prove that the the standard
invariant of the group type subfactors
$M^H \subset M \rtimes K$ has property T if and
only if the group $G$, generated by $H$ and $K$ in $\Out \, M$ has
Kazhdan's property T. This is done without explicitly computing the
symmetric enveloping inclusion associated to $M^H \subset M \rtimes K$.
The key fact here is that property T has certain hereditary properties,
which are established in theorem 2.6. The result regarding
$M^H \subset M \rtimes K$ as mentioned above follows then from
these hereditary properties and the fact that a certain locally trivial
subfactor associated to the kernel $G$ has property T in the sense
of [Po6] if and only if the group $G$ has Kazhdan's property T.
We finish section 2 with a number of explicit examples of 
subfactors with property T standard lattice (see corollary 2.9 and remarks).

\bigskip
\medskip

\heading
1. Property T 
\endheading
\bigskip

We present in this section three related concepts of property T used in this 
paper and fix the notation. Recall first that a countable discrete group $G$ 
has {\it Kazhdan's property T} [K] if, roughly speaking, every unitary 
representation $\pi : G \to \Cal U(\Cal H)$ 
on a Hilbert space $\Cal H$ that has almost invariant vectors has nonzero
invariant vectors. More precisely, $G$ has property T, if
there is a finite set $F \subset G$ and there is
an $\eps > 0$ such that whenever $\pi : G \to \Cal U(\Cal H)$ is a unitary
representation of $G$ on the Hilbert space $\Cal H$ and $\xi \in \Cal H$
is a unit vector with  $\| \pi(s)\xi - \xi \| < \eps $, for all
$s \in F$, then there is a nonzero vector $\eta \in \Cal H$ with
$\pi(g)\eta = \eta$, for all $g \in G$ (see for instance
[HV] for other equivalent definitions).

Next we explain two notions of property T related to inclusions of
von Neumann algebras. Throughout this paper 
$N \subset M$ will denote an extremal inclusion of II$_1$ factors with
finite Jones index unless otherwise stated. As usual,
we denote by $e_N$ the orthogonal projection $L^2(M) \to L^2(N)$
and by $J : L^2(M) \to L^2(M)$ the modular conjugation.
The standard invariant of $N \subset M$, consisting of the
system of higher relative commutants associated to $N \subset M$,
will be denoted by $\Cal G_{N,M}$ (see for instance [Po2]).
We refer to $\Cal G_{N, M}$ also as the {\it standard lattice}
associated to $N \subset M$ (see [Po4]). Inclusions whose index is 
not necessarily finite 
will usually be denoted by $\Cal N \subset \Cal M$.

\medskip
\definition{Definition 1.1}([A-D], [Po1])
Let $\Cal N \subset \Cal M$ be an inclusion of von Neumann algebras and 
let $\Cal M$ be a II$_1$ factor. We say that $\Cal M$ {\it has property T 
relative to} $\Cal N$,
or that the {\it pair} $\Cal N \subset \Cal M$ {\it has property T}, 
if there is an
$\eps > 0$ and operators $x_1$, $x_2, \dots $, $x_n \in \Cal M$ such that
if $\Cal H$ is an $\Cal M$-$\Cal M$ bimodule and $\xi \in \Cal H$ is a vector
satisfying $\| \xi \| = 1$, $\| x_i \xi - \xi x_i \| < \eps$, 
$1 \le i \le n$, $[\xi, y] = 0$, for all $ y \in \Cal N$, then there is
a nonzero vector $\eta \in \Cal H$ such that $[\eta, x ] = 0$ for
all $x \in \Cal M$. $\eta$ is called an $\Cal M$-{\it central vector}.
We call $\{\eps, x_1, \dots , x_n \}$ a {\it critical set} for $\Cal N 
\subset \Cal M$.
\enddefinition
\medskip

Note that we do not assume that  $\Cal N \subset \Cal M$ has finite Jones 
index [Jo]. 
In fact, if $ N \subset M$ is an inclusion of II$_1$ factors with finite 
index, then, using an orthonormal basis
of $M$ over $N$ [PiPo], one sees easily that $M$ has property T relative
to $N$. Let us also mention that if $G$ is a
countable discrete group with a properly outer action on the II$_1$
factor $N$, then the crossed product $N \rtimes G$ has property T relative
to $N$ if and only if the group $G$ has Kazhdan's property T ([A-D], [Po1]).
 
We recall next the notion of property T for the standard invariant
of an extremal subfactor $N \subset M$ as introduced in [Po5], [Po6]. 
If $N \subset M$ is an inclusion of II$_1$ factors with finite Jones
index, we associate to it the {\it symmetric enveloping inclusion}
$M \vee M^{\text{\rm op}} \subset M \bt_{e_N} M^\op$,
where $ M \bt_{e_N} M^\op$ is the II$_1$ factor obtained
in the following way (see [Po3], [Po6] for details):
We consider the C$^*$-algebra $C^*(M,e_N,JMJ) \subset B(L^2(M))$,
generated by $M$, $JMJ$ and $e_N$ on $L^2(M)$.
It can be shown that $C^*(M,e_N,JMJ)$ has a unique trace $\tr$
and one defines 
$M \bt_{e_N} M^\op$ to be the von Neumann algebra obtained via
the GNS construction with respect to $\tr$ from $C^*(M,e_N,JMJ)$.
We can now give the desired definition:

\medskip
\definition{Definition 1.2}[Po6] Let $N \subset M$ be
an extremal inclusion of II$_1$ factors with finite Jones index and with 
standard lattice
(or standard invariant) $\Cal G_{N,M}$. We say that
$\Cal G_{N,M}$ {\it has property T} if the symmetric enveloping
II$_1$ factor $M \bt_{e_N} M^\op$ has property T relative
to the subfactor $M \vee M^{\text{\rm op}}$ (in the sense of definition
1.1). Equivalently, we will say that {\it the standard invariant} 
(or {\it the standard lattice}) of $N \subset M$ {\it has property T}.
\enddefinition
\medskip

It is easy to see that the index of the inclusion 
$M \vee M^{\text{\rm op}} \subset M \bt_{e_N} M^\op$
is finite if and only if $N \subset M$ has finite depth. Thus the
standard lattice of finite depth subfactors has property T as 
expected.

Observe that the definition of property T for $\Cal G_{N,M}$ depends
a priori on the subfactor $N \subset M$. We have however 
the following result [Po6, section 9]: 

\medskip
\proclaim{Theorem 1.3} Suppose $\Cal G_{N,M}$ is the standard lattice
of an extremal subfactor. The property T of 
$\Cal G_{N,M}$ depends only on the (trace preserving) isomorphism
class of $\Cal G_{N,M}$, not on the extremal subfactor from which
it was constructed.
\endproclaim
\medskip

If we let $P$ be a II$_1$ factor and $\theta_0 = id$, $\theta_1, \dots,$
$\theta_n$ automorphism of $P$, then we can consider the locally
trivial inclusion of II$_1$ factors $N^{\theta} \subset M^\theta$,
where $N^\theta = \{ \sum_{i=0}^n \theta_i(x) e_{ii} \, | \, x \in P \, \}$,
$M^\theta = P \otimes M_{n+1}(\Bbb C)$, $\{ e_{ij} \}_{0 \le i, j \le n}$,
the matrix units in $M_{n+1}(\Bbb C)$. These reducible inclusions of
II$_1$ factors provide examples of subfactors with property T standard
invariant [Po6, proposition 9.7]:

\medskip
\proclaim{Proposition 1.4} The standard 
lattice of the subfactor $N^{\theta} \subset M^\theta$ 
has property T if and only if the group $G$ generated by 
the $\theta_i $'s, $0 \le i \le n$, in the outer automorphism
group of $P$ has Kazhdan's property T.
\endproclaim
\medskip

The above proposition is proved by computing the symmetric
enveloping inclusion associated to $N^{\theta} \subset M^\theta$
explicitly. It turns out that it is a crossed product by the group
$G$.

\bigskip

\heading
2. Group type inclusions with property T
\endheading
\medskip

We study in this section the group type inclusions
$M^H \subset M \rtimes K$, where $H$ and $K$ are finite groups
with a properly outer action on the II$_1$ factor $M$ [BH].
$M^H$ denotes as usual the fixed point algebra under the $H$-action
and $ M \rtimes K$ the crossed product.
Note that $M^H \subset M \rtimes K$ has Jones index $|H| \cdot |K|$
and is extremal.
It was shown in [BH] that various analytical properties of these
subfactors, such as finite depth, amenability and strong amenability
(in the sense of [Po2]), are determined by the group
$G = \langle H,K \rangle$, generated by $H$ and $K$ in the outer
automorphism group of $M$. Furthermore, 
$M^H \subset M \rtimes K$ is irreducible if and only if
$H \cap K = \{ e \}$ in $\Out\, M = \Aut \, M / \Int \, M$. 
Recall that if $M$ is the
hyperfinite II$_1$ factor $R$, then any group $G$, which is
a quotient of the free product of $H $ and $K$ gives rise to
an irreducible subfactor of the above type.
We will show in this section the following theorem:

\medskip
\proclaim{Theorem 2.1} Let $H$ and $K$ be finite groups with
properly outer actions $\sigma$, resp. $\rho$ on the II$_1$
factor $M$. The standard invariant of the subfactor 
$M^H \subset M \rtimes K$ has
property T if and only if the group
$G$ generated by $\sigma(H)$ and $\rho(K)$ in the outer autmorphism
group of $M$ has Kazhdan's property T.
\endproclaim
\medskip

This theorem will follow from a more general result, theorem 2.6 below.
We will also present at the end of this section some explicit examples 
of irreducible, infinite depth subfactors whose standard lattice has 
property T.

The next result can be found in [Po1]. We include the proof here
for the convenience of the reader:

\proclaim{Proposition 2.2} 
\roster
\item"i)" Let $\Cal Q \subset \Cal N \subset \Cal M$ be an inclusion 
of II$_1$ factors
with $[\Cal N: \Cal Q] < \infty$. Then $\Cal M$ has property T 
relative to $\Cal Q$ if and
only if $\Cal M$ has property T relative to $\Cal N$.
\item"ii)" Let $\Cal N \subset \Cal M$ be an inclusion of II$_1$ factors 
and let $p \in \Cal N$ be a nonzero projection. Then $\Cal M$ has property 
T relative to $\Cal N$ if
and only if $p\Cal M p$ has property T relative to $p \Cal N p$. 
\endroster
\endproclaim
\demo{Proof} Let us first prove i). If $\Cal M$ has property T relative
to $\Cal Q$, then it is trivial that $\Cal M$ has property T relative 
to $\Cal N$ - the
same critical set works. Conversely, suppose $\Cal M$ has property T relative
to $\Cal N$ and let $\{n_i \}_{i \in I}$ be a (finite) orthonormal basis of
$\Cal N$ over $\Cal Q$. Let $\{\eps, x_1 , \dots , x_n \}$, $x_i \in \Cal M$, 
$1 \le i \le n$,
be a critical set for $\Cal N \subset \Cal M$ (definition 1.1). 
Suppose $\Cal H$ is an $\Cal M$-$\Cal M$-bimodule, $\xi \in \Cal H$, 
$\| \xi \| = 1$ with $[\xi , \Cal Q ] = 0$. Set 
$\eta = \sum_{i \in I} n_i \xi n_i^* \in \Cal H$. Then
we have for all $x \in \Cal N$
$$\align
x \eta &= \sum_i x n_i \xi n_i^* = \sum_j \sum_i n_j 
E_{\Cal Q}(n_j^*xn_i)\xi n_i^*
= \sum_j \sum_i n_j \xi E_{\Cal Q}(n_j^*xn_i)  n_i^* \\
&= \sum_j n_j \xi \big( \sum_i E_{\Cal Q}(n_j^*xn_i)n_i^* \big) 
= \sum_j n_j \xi n_j^* x = \eta x.
\endalign
$$
Furthermore $\| n_i \xi n_i^* - \xi n_i n_i^* \| \le
\| n_i \xi - \xi n_i \| \| n_i^* \| \le [\Cal N:\Cal Q] \| n_i \xi - \xi
n_i \|$, for all $i$. Thus, if we let $\delta = \underset{i} \to \max 
\| n_i \xi - \xi n_i \|$, then

$$\align
\| \eta - [\Cal N: \Cal Q] \xi \| &=
\| \sum_i n_i \xi n_i^* - \sum_i \xi n_i n_i^* \| \le
 \sum_i \| n_i \xi n_i^* - \xi n_i n_i^* \| \\
&\le [\Cal N: \Cal Q] \sum_i \| n_i \xi - \xi n_i \|  \le [\Cal N: \Cal Q]
|I|\delta
\endalign
$$

Hence $\| \eta \| \ge [\Cal N: \Cal Q](1 - |I|\delta)$, so that
$\eta \ne 0$, if $0< \delta < \frac{1}{|I|}$. 
Fix such a $\delta$ and pick an $\eps' > 0$ with
$\eps' < \min \{ \delta, \frac{\eps (1-|I| \delta) }{3 |I|}\}$. Then
$\{ \eps' , x_i, x_i n_j , n_j^*x_i, n_j \}_{ 1 \le i \le n, j \in I }$ is a 
critical set for $\Cal Q \subset  \Cal M$, since we obtain 
$\| x_i \eta ' - \eta ' x_i \| < \eps $, $1 \le i \le n$, 
where $\eta ' = \frac{\eta}{\| \eta \|}$ and
$\{ \eps , x_1, \dots , x_n \}$ was the critical set for 
$\Cal N \subset \Cal M$.

We prove next ii). Let us assume first that $p \Cal M p$ has 
property T relative
to $p \Cal N p$. Let $\{ \eps, y_1, \dots , y_n \}$, $y_i \in p \Cal M p$, 
be a critical set. We may assume $\| y_i \| \le 1$, $1 \le i \le n$.
Since $\Cal M$ is a factor, we can find isometries $v_1, \dots , v_m \in 
\Cal M$
such that $v_1 v_1^* \le p$, $v_i v_i^* = p$, $2 \le i \le m$ and
$\sum_{i=1}^m v_i^*v_i = 1_{\Cal M}$. Let $\delta = \min \{ \eps /6, 
\frac{1}{8m} \}$ 
and set $\{ \delta , x_i \} = \{ y_j \}_{1 \le j \le n} 
\cup \{ v_j \}_{1 \le j \le m} \cup \{ p \}  \subset \Cal M$. 
We will show that this 
consitutes a critical set for
$\Cal N \subset \Cal M$. Let $\Cal H$ be an $\Cal M$-$\Cal M$ bimodule, 
$\xi \in \Cal H$,
$\| \xi \| = 1$ with $[\xi , \Cal N] = 0$ and $\| [\xi , x_i ] \| < \delta$, 
for all $i$. 
Set $\xi_0 = p \xi p $. Clearly, $\xi_0$ is central
for $p \Cal N p$. Note that $\| y_j \xi_0 - \xi_0 y_j \| = 
\| y_j p \xi p - p \xi p y_j \|  \le \| y_j \xi p - y_j \xi \|
+ \| y_j \xi  - \xi y_j \| + \| \xi y_j - p \xi y_j \|
\le 2 \| \xi p - p \xi \| + \| y_j \xi - \xi y_j \| < 3 \delta < \eps$. 

Furthermore 
$\| \xi \|^2 =  \sum_{i=1}^m \| v_i^* v_i \xi \|^2  \le
\sum_i \big( \| v_i^*v_i \xi - v_i^* \xi v_i \|^2 + 2 \| v_i v_i^*\xi
v_i \| \| v_i \xi - \xi v_i \| + \| v_i^* \xi v_i \|^2 \big)$.
Hence $\| \xi \|^2 \le \delta(\delta +2)m + \sum_i \| v_i^* \xi v_i \|^2$.
But $\| v_i^* \xi v_i \|^2 \le \| \xi_0 \|^2$, for all $i$, so that
$\| \xi_0 \|^2 \ge \frac{1}{m} \| \xi \|^2 - \delta(\delta +2)$, which
is $> 0$, if $0 < \delta < \frac{1}{4m}$. 

Since $\| y_j \xi_0 - \xi_0 y_j \| < \eps$, $1 \le j \le n$, there
is a vector $\eta_0 \in p \Cal H p$, $\eta_0 \ne 0$, with
$[\eta_0 , p \Cal M p ] = 0$. Set $\eta = \sum_{i=1}^m v_i^* \eta_0 v_i$,
then it is easy to check that $x \eta = \eta x$, for all $x \in \Cal M$.
Since $\| \eta \|^2 \ge (m-1) \| \eta_0 \|^2$, we have that
$\eta \ne 0$. Thus $\Cal M$ has property T relative to $\Cal N$.

Conversely, suppose that $\Cal M$ has property T relative to $\Cal N$.
Let $\{ \eps , x_1, \dots , x_n \}$, $x_i \in \Cal M$, $1 \le i \le n$,
be a critical set for $\Cal N \subset \Cal M$.
Let $\Cal H_0$ be an $p \Cal M p$-$p \Cal M p$ bimodule and let 
$\xi_0 \in \Cal H_0$,
$\| \xi_0 \| =1$, with $[\xi_0 , p \Cal N p]=0$. As above, we can find 
partial
isometries $v_1, \dots , v_m \in \Cal M$ with $v_i v_i^* =p$, $2 \le i \le m$,
$v_1 v_1^* \le p$ and $\sum_{i=1}^m v_i^* v_i = 1$. Consider the
set $\{ \delta, v_r x_k v_s^* \}_{1 \le k \le n, 1 \le r,s \le m}$ -
it will be a critical set for $p \Cal N p \subset p \Cal M p$, 
if $\delta $ is
small enough. Set $\alpha =  \tr(p)^{-1}$ and let 
$\Cal H$ be the $\alpha$-amplification of $\Cal H_0$. The vector
$\xi = \sum_{i=1}^m v_i^* \xi v_i \in \Cal H$ is $\Cal N$-central and
$\ne 0$. We leave it as an exercise to check that
$\| \xi x_k - x_k \xi \| < \eps$, $ 1 \le k \le n$, if $\delta $
is chosen small enough. Thus, since $\Cal M$ has property T relative
to $\Cal N$, there is an $\Cal M$-central vector $\eta \in \Cal H$, 
$\eta \ne 0$.
But then $p \eta p \ne 0$ and $p \eta p \in p \Cal H p = \Cal H_0$ is
the desired $p \Cal M p$-central vector in $\Cal H_0$. \qed 
\enddemo

Proposition 2.2 implies the following corollary (see [Po6] for a 
slightly more general statement):

\proclaim{Corollary 2.3}
Let $N \subset M$ be an extremal inclusion 
of II$_1$ factors with finite index
and let $N \subset M \overset{e_1}\to\subset M_1$ be the basic 
construction. Then the standard invariant of $N \subset M$ has property T 
if and only if the standard invariant of $M \subset M_1 $ has property T.
\endproclaim
\demo{Proof}
Consider the inclusions 
$M \vee N^{\text{\rm op}} \subset M_1 \vee  N^{\text{\rm op}} \subset
M \bt_{e_N} M^\op$ and
$M \vee N^{\text{\rm op}} \subset M \vee  M^{\text{\rm op}} \subset
M \bt_{e_N} M^\op$. By proposition 2.2, i) we have that
$M \bt_{e_N} M^\op$ has property T relative to $M \vee N^{\text{\rm op}}$
if and only if it has property T relative to $M_1 \vee  N^{\text{\rm op}}$.
Using the second inclusion, we get similarly that
$M \bt_{e_N} M^\op$ has property T relative to
$M \vee N^{\text{\rm op}}$ if and only if it has property T relative
to $M \vee  M^{\text{\rm op}}$. Thus $M \bt_{e_N} M^\op$ has
property T relative to  $M \vee  M^{\text{\rm op}}$ (i.e.
$\Cal G_{N,M}$ has property T) if and only if 
$M \bt_{e_N} M^\op$ has property T relative to $M_1 \vee  N^{\text{\rm op}}$.
But $M_1 \vee  N^{\text{\rm op}} \subset M \bt_{e_N} M^\op$
is a reduced inclusion of the symmetric enveloping inclusion
for $M \subset M_1$, reduced by $e_1^{\text{\rm op}}$ [Po6, section 2].
By proposition 2.2 ii) we are done.  \qed
\enddemo

We also need the following observation from [Po6]:

\bigskip 
\proclaim{Lemma 2.4} Let $Q \subset N \subset M$ be an extremal inclusion
of II$_1$ factors with finite index and let $\{ m_j \}_{j \in J}$ be an 
orthonormal basis
of $N$ over $Q$. The map which sends $e_N$ to
$\sum_j m_j e_Q m_j^*$ and which is the identity on
$\text{\rm Alg}(M,M^{\op})$ implements a unital embedding of
$C^*(M,e_N,JMJ)$ into $C^*(M,e_Q,JMJ)$ and thus a unital embedding of
$M \bt_{e_N} M^\op$ into $M \bt_{e_Q} M^\op$.
\endproclaim

\demo{Proof}
If $x \in N$, then we have $x = \sum_j m_j E_Q(m_j^*x)$ and hence
$\sum_j m_j e_Q m_j^* (\hat x)$ $ = \sum m_j E_Q(m_j^*x)^\wedge$ $= \hat x$.
Also, if $x \in M$, $\hat x \perp N$, then $\widehat{m_j^* x} \perp N$
as well, so $\sum m_j e_Q m_j^* (\hat x) = \sum_j m_j e_Q (\widehat
{m_j^* x})$ $=\sum_j m_j e_Q e_N (\widehat{m_j^* x}) = 0$.
Hence $\sum_j m_j e_Q m_j^* = e_N$ in $B(L^2(M))$. 
Thus we get
$C^*(M,e_N,JMJ) \subset C^*(M,e_Q,JMJ)$. By [Po6]
it follows that this implements an inclusion
$M \bt_{e_N} M^\op \subset M \bt_{e_Q} M^\op$ as well.
\qed
\enddemo
\medskip 

\proclaim{Lemma 2.5} Let $N \subset P \subset M$ be an extremal inclusion
of II$_1$ factors with finite index 
and let $P_{-1} \subset N \subset P \subset M$ be such that $N \subset
P$ is the basic construction for $P_{-1} \subset N$. 
Then $\Cal G_{N,M}$ has property T if and only if $\Cal G_{P_{-1},M}$
has property T.
\endproclaim

\demo{Proof} By lemma 2.4 we have the inclusion
$(M \vee M^\op \subset M \bt_{e_N} M^\op)$ $\subset
(M \vee M^\op \subset M \bt_{e_{P_{-1}}} M^\op)$.
We will show that it is in fact an equality. 
Clearly, the statement follows then by applying proposition 2.2 i).
Let $N_1 \subset P_{-1} \subset N \subset P \subset M$
be such that $N \subset M$ is the basic construction for
$N_1 \subset N$ and $N \subset P$ is the basic construction
for $P_{-1} \subset N$ (see for instance [Bi]). Then by 
lemma 2.4 we
have the inclusions $C^*(M,e_N,J_MMJ_M) \subset C^*(M,e_{P_{-1}},J_M M J_M)
\subset C^*(M,e_{N_1},J_M M J_M)$. By [Po6, section 2] the first algebra
actually equals the last (using the fact that the Jones projection
$M \to N_1$ can be written in terms of the $e_i$'s).  \qed
\enddemo

We show next the main theorem of this section.

\proclaim{Theorem 2.6} Let $N \subset P \subset M$ be an extremal
inclusion of II$_1$ factors with finite index and let 
$P_{-1} \subset N \subset P \subset M \subset P_1$ be such that 
$N \subset P$ is a basic construction
for $P_{-1} \subset N$ and $M \subset P_1$ is a basic construction
for $P \subset M$. Then the standard invariant of $P_{-1} \subset P_1$ 
has property T if and only if the standard invariant of $N \subset M$ 
has property T.
The same statement holds if ``property T'' is replaced by ``amenable''.
\endproclaim

\demo{Proof} Let 
$P_{-2} \subset N_1 \subset P_{-1} \subset N \subset P \subset M
\subset P_1$ be such that 
$P_{-2} \subset N \subset P_1$ and $N_1 \subset N \subset M$
are basic constructions. Thus the standard invariant of 
$P_{-2} \subset N$ has property T if and only if the standard invariant
of $N \subset P_1$ has property T (corollary 2.3). Furthermore,
by lemma 2.5 $\Cal G_{N,M}$ has property T if and only if 
$\Cal G_{P_{-1},M}$ has property T. By corollary 2.3 we have
that $\Cal G_{N,M}$ has property T if and only if
$\Cal G_{N_1,N}$ has property T and the latter has property T
if and only if $\Cal G_{P_{-2},N}$ does (lemma 2.5 again). 
Hence, $\Cal G_{N,M}$ has property T if and only if
$\Cal G_{N,P_1}$ has property T.

Now consider the inclusions $P_{-1} \subset N \subset P \subset P_1$
and apply lemma 2.5 again. Thus $\Cal G_{P_{-1},P_1} $ has
property T if and only if $\Cal G_{N,P_1}$ has property T. Together
with the above equivalences we have therefore  that
the standard invariant of $P_{-1} \subset P_1$ has property T
if and only if the one of $N \subset M$ has property T.  

Using ([Po2], [Po6]) the same proof goes through if
``property T'' of the standard invariant is replaced by ``amenability''
of the standard invariant
(assuming of course in addition that the II$_1$ factors appearing are 
hyperfinite, if one replaces ``property T'' by  ``amenable inclusions'' 
rather than ``amenable standard invariant''). \qed
\enddemo
\medskip

To prove theorem 2.1 we also need the following proposition.

\proclaim{Proposition 2.7} Let $H$ and $K$ be finite groups with
properly outer actions $\sigma$, resp. $\rho$ on the II$_1$ factor $P$.
Denote $N = P^H$ and $M = P \rtimes K$ and let
$P_{-1} \subset N \subset P \subset M \subset P_1$ be such that
$N \subset P$ is the basic construction for $P_{-1} \subset N$
and $M \subset P_1 $ is the basic construction for $P \subset M$.
Then $P_{-1} \subset P_1$ is a locally trivial inclusion given
by $\{ \rho_k\sigma_h \}_{h \in H, k \in K}$.
\endproclaim

\demo{Proof}
Note that $P \subset P_1$ can be described as follows: $P_1
= P \otimes B(L^2(K))$ with matrix units $\{ e_{k,k'} \}_{k,k' \in K}$
and the inclusion $P \hookrightarrow P_1$ is given by
$P \ni x \to \sum_{k \in K} \rho_k(x)e_{k,k}$.
By perturbing each $\sigma_h$ and $\rho_k$ by some inner
automorphism of $P$ we can assume that they all leave the same
$B \subset P$ invariant, with matrix factor
$B \cong B(l^2(H)) $ and matrix units 
$\{ f_{h,h'} \}_{h,h' \in H}$. But
then $P_{-1} \subset P$ is isomorphic to the inclusion
$P_{-1} \hookrightarrow P_{-1} \otimes B(l^2(H))$, given by
$P_{-1} \ni x \to \sum_{h \in H} \sigma_h(x) f_{h,h}$.
Altogether this gives that $P_{-1} \subset P_1$ is isomorphic
to the inclusion $P_{-1} \hookrightarrow P_{-1} \otimes
B(l^2(H)) \otimes B(l^2(K))$ given by
$P_{-1} \ni x \to \underset{h \in H, k \in K}\to \sum \rho_k\sigma_h(x)
f_{h,h}e_{k,k}$.  \qed
\enddemo

Theorem 2.1 follows now immediately from proposition 1.4 and theorem 2.6.

\medskip

\remark{Remark 2.8} Observe that we recovered in theorem 2.6 and
proposition 2.7 the result
[BH, theorem 4.9]) showing that $P^H \subset P \rtimes K$ has amenable 
principal graph if and only if the group $G$ generated by $\sigma(H)$ and
$\rho(K)$ in $\Out (P)$ is amenable.
\endremark

\bigskip

Numerous explicit examples of infinite groups $G$ that have property T and 
are quotients of the free product of two finite groups $H$ and $K$ can be 
found in [BS] (we would
like to thank Pierre de la Harpe for pointing out this reference to us).
The construction of these groups is based on Sarnak's
examples of Ramanujan graphs [Sa, chapter 3] and an explicit 
presentation of this class of property T groups can be found in
[BS, theorem 2 and corollary 2]. Since any discrete group has a properly outer
action on the hyperfinite II$_1$ factor $R$, we obtain therefore (by applying
theorem 2.1) irreducible, infinite depth subfactors of the form 
$R^H \subset R \rtimes K$ whose standard invariant has property T.

We will construct next irreducible, infinite depth subfactors 
with Jones index $6$ whose standard invariant has property T. By theorem 2.1 we 
need to find an infinite property T group which is a quotient of
$PSL(2,\Bbb Z) = \Bbb Z_2 * \Bbb Z_3$. 

\proclaim{Corollary 2.9} There are irreducible, infinite depth
subfactors of the hyperfinite II$_1$ factor with Jones index $6$
whose standard lattice has property T.
\endproclaim

\demo{Proof} Let $G$ be the group given by the presentation
$\langle x , y \, | \, x^2 = y^3 = (xy)^{12} = (xy^{-1}xy^{-1}xyxy
xy^{-1}xy)^2=1 \rangle$. It is shown in [Co] that this group is
a finite index subgroup of $SL(3,\Bbb Z) \rtimes Z_2$, where $Z_2$
acts on $SL(3,\Bbb Z)$ as the inverse-transpose automorphism.
Since $SL(3,\Bbb Z)$ has property T ([K], see also [HV]) and property
T is preserved under finite extensions and finite index subgroups
(see for instance [HV]), the group $G$ must have property T as well.
$G$ is the desired infinite quotient of $PSL(2,\Bbb Z)$ with
property T - apply then theorem 2.1.  \qed 
\enddemo

\remark{Remark 2.10}\roster 
\item"i)" In [CRW] $SL(3,\Bbb Z)$ is written as a quotient of
$A_4 * \Bbb Z_3$. By theorem 2.1 we obtain therefore an irreducible,
infinite depth subfactor with index $36$ and property T standard 
invariant.

\item"ii)"In [CRW] $SL(3,\Bbb Z) \rtimes \Bbb Z_2$ is written as a 
quotient of $(\Bbb Z_2 \times \Bbb Z_2 ) * \Bbb Z_3$. Again by
theorem 2.1, we obtain an irreducible, 
infinite depth subfactor with index $12$ and property T standard 
invariant and intermediate subfactors of index $2$ and $6$.
\endroster
\endremark

We end this paper by mentioning a few problems that are closely related to
the above work and seem to be of interest at this point.
\roster
\item "a)" Are there irreducible, infinite
depth subfactors with property T standard lattice and Jones index $< 6$? 
Observe that such a subfactor would have necessarily index
$> 3 + \sqrt{3}$ by [Ha] and [Po6, section 9].

\item "b)" The group type inclusions $M^H \subset M \rtimes K$ have
intermediate subfactors. Construct irreducible, 
infinite depth subfactors whose standard lattice has property T that 
do {\it not} have intermediate subfactors. 

\item "c)" Construct irreducible, infinite
depth subfactors with {\it noninteger} Jones index whose standard 
lattice has property T. 

\item "d)" It would be interesting to compute the symmetric enveloping inclusion
associated to $M^H \subset M \rtimes K$ and to give then a different
proof of our theorem 2.1.

\item "e)" Suppose $N \subset P \subset M$ is an extremal inclusion of
II$_1$ factors with finite index such that the standard lattice of 
$N \subset M$ has property T. Suppose furthermore that
one of the intermediate inclusions has finite depth. Does it then follow
that the standard lattice of the other intermediate inclusion has   
property T ? The answer may very well be negative.

\endroster

\vfill
\eject

\enddocument
 
\vfill
\bye